# SOME CLASSES OF DISTRIBUTIONS ON THE NON-NEGATIVE LATTICE

**S. Satheesh and N. Unnikrishnan Nair**
Cochin University of Science and Technology

**Abstract**

A method for constructing distributions on the non-negative lattice of points $I_o = \{0,1,2, ....\}$ as discrete analogue of continuous distributions on $[0,\infty)$ is presented. A justification of the definition of discrete class-L laws is provided. Discrete analogue of distributions of the same type and the role of Bernoulli law in this context is discussed. Generalizations of some distributions and properties of $\alpha$-Poisson laws are given. The geometric compounding problem for discrete distributions is studied by introducing discrete semi Mittag-Leffler laws.

## 1 Introduction

The discussions in literature on the properties of distributions using characteristic functions (CF) or Laplace transforms (LT) do not clearly specify whether such a framework can be made use of in studying discrete distributions as well. For example, it is not explicitly established whether the concepts of geometric compounding and distributions of the same type naturally carry over to the discrete domain. Some important work done in literature on continuous distributions, which opens up the possibility of discussion on discrete distributions along these lines are: Kakosyan, et. al. (1984), Lin (1994) and Sandhya (1991a,b) on geometric compounding. Certain studies of discrete analogue of continuous distributions using probability generating functions (PGF) are by Steutel and van Harn (1979) on discrete class-L laws and by Jayakumar and Pillai (1992) on discrete Mittag-Leffler laws.

In the present paper our objective is to construct PGFs of certain discrete distributions on the non-negative lattice $I_o = \{0,1,2, ....\}$, as discrete analogue of the LTs of continuous distributions on $[0,\infty)$ and discuss distributions of the same type and geometric compounding in the discrete domain.

In Section 2 we present a characterization result that relates PGFs to LTs, which enable us to construct a variety of PGFs as discrete analogues of their LT counter parts. We then show that a formal analogue of class-L laws lead to the definition of discrete class-L laws by Steutel and van Harn (1979). In Section 3 the concept of discrete analogue of distributions of the same type is discussed and the role of Bernoulli laws in this context is highlighted. Following this, generalizations of Bernoulli and Poisson laws and same known results are presented in Section 4 along with some properties of $a$-Poisson laws. Finally, in Section 5 we define and study discrete semi Mittag-Leffler laws and its subclass, the discrete Mittag-Leffler laws in various compound geometric setups.

## 2 Basic Results

**Lemma 2.1** : If $\phi(s)$ is a LT, then $P(s) = \phi(1-s)$, $0<s<1$ is a PGF. Conversely, if $P(s)$ is a PGF and $P(1-s)$ is completely monotone for all $s>0$, then $\phi(s) = P(1-s)$ is a LT.

**Proof** : We have $\phi(0) = 1 \Leftrightarrow P(1) = 1$ and thus the norming requirements for LTs and PGFs are satisfied by the construction. Again, since $\phi(s)$ is a LT it is completely monotone (CM) for all $s>0$ and hence $P(s) = \phi(1-s)$ is absolutely monotone (AM) for all $0<s<1$. Thus by Feller (1966, p.221) $P(s)$ is a PGF. Conversely, if $P(s)$ is a PGF it is AM for all $0<s<1$ and hence $\phi(s) = P(1-s)$ is CM for all $0<s<1$. But for $\phi(s)$ to be a LT it must be CM for all $s>0$ (Feller (1966), p.415). Hence the proof is complete.

The following example shows that $P(1-s)$ need not be CM for all $s>0$.

**Example 2.1** : Starting from the PGF of the Bernoulli law, $[1-b(1-s)]$, $0<b<1$ it follows that $[1-b(1-s)^\alpha]$, $0<\alpha \leq 1$ is also a PGF which we refer to as that of an $\alpha$-Bernoulli law. Setting

$$P(s) = 1 - b(1-s)^\alpha, \ 0<s<1, \ 0<\alpha \le 1, \ 0<b<1.$$

$P(1-s) = 1 - bs^\alpha$ which is CM for $0<s<1$. But when $s>1$, $P(1-s)$ could be negative. Thus $P(1-s)$ is not CM for all $s>0$. Hence $P(1-s)$ is not a LT.

Next, we give a formal justification of discrete class-L laws.

**Theorem 2.1** : A PGF P(s) is in discrete class-L if and only if, for each $0<\alpha<1$, there exists another PGF $P_\alpha(s)$ such that

$$P(s) = P(1-\alpha+\alpha s) \, P_\alpha(s). \tag{1}$$

**Proof** : Class-L laws on $[0,\infty)$ are defined by LTs $\phi(s)$, satisfying, for each $0<\alpha<1$,

$$\phi(s) = \phi(\alpha s) \, \phi_\alpha(s) \tag{2}$$

where $\phi_\alpha(s)$ is another LT. In terms of PGFs (constructed by Lemma.2.1 from these LTs) equation (2.2) reads, for each $0<\alpha<1$,

$$P(1-s) = P(1-\alpha s) \, P_\alpha(1-s), \ 0<s<1 \tag{3}$$

Setting $1-s = u$ in (2.3) we get (2.1).

**Corollary 2.2** : If a LT $\phi(s)$ is in class-L then the PGF $P(s) = \phi(1-s)$ is in discrete class-L. Conversely if $P(s)$ is in discrete class-L and $\phi(s) = P(1-s)$ is a LT, then $\phi(s)$ is in class-L.

Discrete analogue of stable laws with LT, $\phi(s) = \exp\{-\lambda s^\alpha\}$, $\lambda>0$, $0<\alpha\le1$ is given by the PGF

$$P(s) = \exp\{-\lambda(1-s)^\alpha\}. \tag{4}$$

When $\alpha = 1$ the Poisson law results as the discrete analogue. Thus (2.4) can be considered as a generalization of Poisson laws, which we refer to as $\alpha$-Poisson. Further properties of this distribution will be discussed in Section 4.

Setting $\alpha = 1$ and randomizing $\lambda$ in (2.4) with a distribution having LT $\phi$ we have $P(s) = \phi(1-s)$ and thus:

**Theorem 2.2** : Every PGF $P(s) = \phi(1-s)$, where $\phi$ is a LT, is a mixture of Poisson laws.

**Note** : For the PGF $P(s) = \phi(1-s)$, where $\phi(s)$ is a LT, the mixing distribution is the one with LT, $\phi(s) = P(1-s)$. Thus every geometric law on $I_o$ is an exponential mixture of Poisson laws. The distribution in Example.2.1 is not a mixture of Poisson laws.

## 3  Distributions of the same D-type

It is known that distributions of the same type in the continuous case (that is, r.v's $X$ and $Y$ satisfying $X = cY$ or their LTs satisfying $\phi_X(s) = \phi_Y(cs)$, for all $s>0$ and some $c>0$) do not have an analogy in the lattice case. Here we arrive at it in terms of PGFs using the construction in Lemma.2.1 and then specialize it to the case $0<c<1$ to define and study the concept of D-type.

**Definition 3.1** : Two PGFs $P_1(s) = \phi_1(1-s)$ and $P_2(s) = \phi_2(1-s)$ (where $\phi_1$ and $\phi_2$ are LTs) are of the same type if and only if $\phi_1(1-s) = \phi_2(c(1-s))$, for all $0<s<1$ and some $c>0$.

Clearly, this definition applies to PGFs derived from LTs, while Example.2.1 (including the Bernoulli law) suggest that this is not the case always. Further, for these two distributions the range $c<1$ alone is safely applicable in Definition.3.1. Exploring this range for $c$ showed that it has some nice implications as is seen below which motivates us to coin the nomenclature D-type in the next definition. Also, this is the range of $c$ in the context of distribution of the same types in summation schemes discussed in Sections 4 and 5 of this paper. However, Definition.3.1 is still relevant which will be highlighted after Example.3.1 as a note.

**Definition 3.2** : Two PGFs $P_1(s)$ and $P_2(s)$ are of the same D-type if and only if $P_1(1-s) = P_2(1-cs)$, for all $0<s<1$ and some $0<c<1$.

**Theorem 3.1** : Two PGFs $P_1(s)$ and $P_2(s)$ are of the same D-type if and only if $P_1$ is a $P_2$ compounded Bernoulli law.

**Proof** : The assertion is proved by setting $u = 1-s$ in Definition.3.2 to get

$$P_1(u) = P_2(1- c+cu), \quad 0<c<1 \tag{1}$$

the Bernoulli probability being $c$.

Notice that if $P_1$ and $P_2$ are constructed from LTs $\phi_1$ and $\phi_2$ as in Lemma.2.1 then writing (3.1) using the structure of the corresponding LTs we have,

$$\phi_1(1-u) = \phi_2[1-(1-c+cu)] = \phi_2(c(1-u)),$$

as it should be in the light of Definition.3.1. From equation (3.1) we further have:

**Theorem 3.2** : Two non-negative lattice r.vs $X$ and $Y$ will have the same D-type distribution if and only if $X = \sum_{i=1}^{Y} Z_i$ for some i.i.d Bernoulli r.vs $\{Z_i\}$ independent of $Y$, the Bernoulli probability being $c$.

The above two theorems also justify the replacement of $cX$ in the continuous set up by $c \circ X = \sum_{i=1}^{Y} Z_i$ to obtain the corresponding lattice analogue, as done in Steutel and van Harn (1979). Next theorem brings out another role of Bernoulli laws.

**Theorem 3.3** : Every PGF $P_1(s) = \phi_1(1-s)$, where $\phi_1$ is a LT, is a compound of Bernoulli laws.

**Proof** : We have

$$P_1(s) = \phi_1(1-s),\ 0<s<1,$$

$$= \phi_1[ab(1-s)],\ 0<b<1,\ ab = 1$$

$$= \phi_1\{a[1-(1-b+bs)]\}$$

$$= P_2(1-b+bs),$$

where $P_2(s) = \phi_1[a(1-s)]$. This completes the proof.

Clearly $P_1$ and $P_2$ here are of the same type according to the Definition.3.1. Also the description of the PGF, $P_2(s) = \phi_1[a(1-s)]$, for $a>0$ and possibly $a>1$ above is justified because $\phi_1$ is a LT. As an illustration consider the following example.

**Example.3.1** : Let $Q_1(s) = 1/\{1+ \lambda(1-s)^\alpha\}$, $\lambda > 0$, $0 < \alpha \leq 1$. This PGF can be derived from the Mittag-Leffler law with LT $1/(1+\lambda s^\alpha)$ by invoking Lemma.2.1. Now choose $a$ and $0 < b < 1$ such that $ab = \lambda$. Then,

$$Q_1(s) = 1/\{1+ ab(1-s)^\alpha\} = \phi_1\{ab(1-s)\}$$

$$= 1/\{1+ a[1-(1-b^{1/\alpha} + b^{1/\alpha}s)]^\alpha\}$$

$$= Q_2(1- b^{1/\alpha} + b^{1/\alpha}s),$$

where $Q_2(s) = 1/\{1+ a(1-s)^\alpha\}$

$$= \phi_1\{a(1-s)\}.$$

**Note** : This is a situation where Definition.3.1 is still relevant as it takes care of the full range of the parameter $\lambda$. Similar is the case with the PGF (2.4) of $\alpha$-Poisson laws.

The following example demonstrates that the conclusion of Theorem.3.3 can still hold good even when Lemma.2.1 does not derive the PGF from a LT.

**Example 3.2** : Let $P_1(s) = 1-\lambda(1-s)^\nu$, $0 < \lambda < 1$ and $0 < \nu \leq 1$ (same as Example.2.1).

Choose $b$ such that $0 < \lambda < b < 1$ and write $ab = \lambda$ so that $0 < a < 1$ also holds true. Now,

$$P_1(s) = 1 - ab(1-s)^\nu$$

$$= 1 - b[a^{1/\nu}(1-s)]^\nu$$

$$= 1 - b[1 - (1 - a^{1/\nu} + a^{1/\nu} s)]^\nu$$

$$= P_2(1- a^{1/\nu} + a^{1/\nu} s), \text{ where}$$

$$P_2(s) = 1 - b(1-s)^\nu.$$

**Remark 3.1** : Another point stressed here is that when $P_1(s) = \phi_1(1-s)$, where $\phi_1$ is a LT, the choice $0 < b < 1$ alone is to be assured and the value of $a$ being greater than or less than unity is immaterial. But in the case of PGFs not constructed from LTs by Lemma.2.1 one should take care that both the factors $a$ and $b$ are less than unity. Eg. in Example.3.2 if $b<\lambda$, then $\lambda = ab$ would imply that $a >1$ and hence speaking about a Bernoulli probability of $a^{1/\nu}$ is meaningless.

## 4 Generalizations of some discrete laws

### 4.1 Generalizations of binomial and Poisson Laws

From the PGF of the $\alpha$-Bernoulli law in Example.2.1 it follows that

$P_n(s) = [1-p(1-s)^\alpha]^n$, $0<p<1$, $0<\alpha \leq 1$, $n = 1, 2, ...$

is another PGF which we refer to as that of the $\alpha$-binomial law. Its limiting case, as: $n\to\infty$, $p\to 0$ such that $np = \lambda$, a constant is $P(s) = \exp\{-\lambda(1-s)^\alpha\}$ which is the PGF of $\alpha$-Poisson law.

As a generalization of geometric laws Jayakumar and Pillai (1992) considered the distribution with PGF, $P(s) = [1+\lambda(1-s)^\alpha]^{-1}$, $\lambda>0$, $0<\alpha\leq 1$ which they called the discrete Mittag-Leffler law. (We have derived this PGF invoking Lemma.2.1 in Example.3.1). Clearly $P_n(s) = \left\{1+\frac{\lambda}{n}(1-s)^\alpha\right\}^{-n}$ is also a PGF. Since

$$\lim_{n\to\infty} P_n(s) = \exp\{-\lambda(1-s)^\alpha\},$$ we have proved:

**Theorem 4.1** : Corresponding to each $\alpha$-Poisson law we have a sequence of discrete Mittag-Leffler laws converging weakly to it.

### 4.2 Discrete Semi Stable and $\alpha$-Poisson Laws

The class of continuous functions

$\psi(s) = a\psi(bs)$, for all $s\in R$ and some $a,b>0$ \hfill (1)

with $\psi(o) = 0$ have been discussed in the context of CFs of semi stable laws and regression equations in Kagan, et. al ((1973), p. 9, 163, 323, 324)), integrated Cauchy functional equations in Pillai and Anil (1996), a variation in approximating the gamma function in Dubec (1990) and in a study of the near-constancy phenomena in branching processes in Biggins and Bingham (1991). It has been proved that for (4.1) to hold, the condition $0<b<1<a$ is necessary and that there must exist a unique $\alpha>0$ such that $ab^{\alpha}=1$. $\alpha$ has to be restricted to $(0,2]$ when $\phi(s) = \exp\{-\psi(s)\}$ is a CF (that of a semi stable law) and $\alpha \in (0,1]$ when $\phi(s)$ is a LT. Pillai and Anil (1996) have shown that the general solution of (4.1) is $\psi(s) = s^{\alpha} h(s)$ where $h(s)$ is a periodic function in $\log(s)$ with period $[-\log(b)]$. Also $\psi(s) = s^{\alpha}$ is clearly such a function.

**Remark 4.1**: In this context two problems are relevant. (i) Whether there are solutions to $\psi(s)$ other than $\psi(s) = s^{\alpha}$, so that functions of $\psi(s)$ are LTs (in this case $\psi(s)$ should have CM derivative), and (ii) conditions under which they reduce to $s^{\alpha}$. Regarding (i), CM solutions to $\xi(as) = b\xi(s)$ were given by Dubec (1990), Biggins and Bingham (1991) and their results clearly show that there are functions $\xi(s)$ other than $s^{-\alpha}$ and under certain conditions $\xi(s)$ is close to a multiple of $s^{-\alpha}$. Bunge (1996) has used these results to describe LTs of semi stable laws as $\exp\{-\psi(s)\}$ where $\psi'(s) = \xi(s)$ and those of semi Mittag-Leffler laws as $1/[1+\psi(s)] = \xi(s)/[1+\xi(s)]$. Also, Jayakumar and Pillai (1993) have given the example of

$$\psi(s) = s^{\alpha}\{1 - A\cos[k\ \log(s)]\}, k = -2\pi/\log(b), 0<b<1 \text{ and } 0<A<1.$$

One can see that $1/[1+\psi(s)]$ is a LT by invoking the one-to-one correspondence between a real CF and a LT in Schoenberg (1938) and noticing that Pillai (1985) had shown that

$$1/[1+ |t|^{\alpha}\{1 - A\cos[k\ \log(|t|)]\}], k = -2\pi/\log(b)$$

is a Polya type CF. Thus definitions involving $\psi(s)$ lead to cases other than $\psi(s) = s^{\alpha}$ so that (ii), the conditions under which they reduce to $s^{\alpha}$ is also important since the corresponding LTs have nice summation stability properties.

Thus we are in a position to develop discrete analogues of semi stable and semi Mittag-Leffler laws by invoking Lemma.2.1 and discuss schemes under which they become discrete stable and discrete Mittag-Leffler (cases corresponding to $s^{\alpha}$).

**Definition 4.1**: In analogy with the continuous case, a PGF P(s) is said to be discrete semi stable (DSS($a,b,\alpha$)), if $P(s) = \exp\{-\psi(1-s)\}$, where

$$\psi(1-s) = a\psi[b(1-s)] \tag{2}$$

for all $0<s<1$ and some $0<b<1<a$ satisfying $ab^\alpha=1$ for a unique $0<\alpha\leq 1$.

**Remark 4.2** : From the inequality $0< \alpha = \dfrac{-\log(1/a)}{-\log b} \leq 1$, a solution for $0<\alpha\leq 1$ of $ab^\alpha = 1$ exist if and only if $ab\leq 1$. This imposes certain restrictions in the choice of Bernoulli probabilities $b$, as will be seen in subsequent deliberations.

**Theorem 4.2** : The sum of $n$ i.i.d discrete variables is distributed as the same D-type as the components, if and only if it is DSS($n,b,\alpha$).

**Proof** : We have

$$P(s) = \exp\{-\psi(1-s)\}$$

$$= \exp\{-n\psi[b(1-s)]\}$$

$$= [\exp\{-\psi[b(1-s)]\}]^n.$$

This proves the if part. The only if part follows by setting $-\log P(s) = \psi(1-s)$ and retracing the above steps. Notice that $P(s) = \phi(1-s)$ do not vanish in the domain $0<s<1$.

**Theorem 4.3** : If a discrete r.v can be expressed as the sum of $n_1$ and $n_2$ independent variables of the same D-type such that $\log n_1/\log n_2$ is irrational, then the variable is $\alpha$-Poisson.

**Proof** : When the DSS($n,b,\alpha$) law can assume, for the same $0<\alpha\leq 1$ two different values for $n$, say $n_1$ and $n_2$ such that their logarithms are in irrational ratio, then $\psi(1-s) = \lambda(1-s)^\alpha$, for some $\lambda>0$ constant (Kagan, et. al (1973),p.9, 323, 324 ) and hence the variable is $\alpha$-Poisson.

**Theorem 4.4** : Let $\{Y_j\}$ be a sequence of i.i.d Bernoulli variables with parameter $b$. Let $M$ be a non-negative lattice variable independent of $\{Y_j\}$ with PGF P($s$) satisfying

$$\lim_{s \to 1} \frac{-\log P(s)}{(1-s)^{\alpha}} = \lambda, \text{ a positive constant and } 0 < \alpha \leq 1. \tag{3}$$

Define $X = \sum_{j=1}^{M} Y_j$. Let $X_1$ and $X_2$ are independent copies of $X$. Then $M$ is identically distributed as $X_1 + X_2$ if and only if $M$ is $\alpha$-Poisson and $b^{\alpha} = \frac{1}{2}$.

**Proof**: Let $P(s) = \phi(1-s)$.

We have $\phi(1-s) = [\phi(b(1-s))]^2$.

Setting $\psi(1-s) = -\log \phi(1-s)$,

$$\exp\{-\psi(1-s)\} = \exp\{-2\,\psi[b(1-s)]\}$$

Solving we find $b^{\alpha} = \frac{1}{2}$.

Writing $Q(1-s) = \dfrac{\psi(1-s)}{(1-s)^{\alpha}}$,

$$\exp\{-(1-s)^{\alpha} Q(1-s)\} = \exp\{-2 \cdot \tfrac{1}{2}(1-s)^{\alpha} Q[b(1-s)]\},$$

or

$$Q(1-s) = Q[b(1-s)], \quad 0 < b < 1$$

$$= Q[b^n (1-s)], \text{ for } n \geq 1 \text{ integer}.$$

Hence by virtue of the condition (4.3) satisfied by $P(s)$, $Q(1-s) = \lambda$ and hence $P(s) = \exp\{-\lambda(1-s)^{\alpha}\}$, and the proof is complete.

**Note**: **(i)** The above result can be extended to the case when $M$ is identically distributed as $\sum_{i=1}^{n} X_i$, in which case $b^{\alpha} = \frac{1}{n}$.

**(ii)** Choosing $Y_j$ to be Bernoulli with parameter $\frac{1}{2}$, in Theorem 4.4 we have a characterization of Poisson ($\lambda$) variable and in this case $X$ also is Poisson by virtue of Raikov's theorem.

**(iii)** The restriction in the choice of the Bernoulli probability is $b \leq \frac{1}{n}$ or ½ as the case may be (c.f Remark.4.1).

**Theorem 4.5** : Under the assumptions of Theorem 4.4, the condition $\sum_{i=1}^{m} M_i \stackrel{d}{=} \sum_{i=1}^{n} X_i$, where $M_i$'s are independent copies of $M$, characterizes the $\alpha$-Poisson law. In this case, the condition $n>m$ should be satisfied.

**Proof** : The condition implies

$$\exp\{-m\,\psi(1-s)\} = \exp\{-n\psi[b(1-s)]\}$$

$$\Rightarrow \quad m\,\psi(1-s) = n\,\psi[b(1-s)]$$

or $\quad \psi(1-s) = \frac{n}{m}\,\psi[b(1-s)]$.

since $\frac{n}{m} > 1$ must be satisfied we should have $n>m$. Proof of $\psi(1-s) = \lambda(1-s)^{\alpha}$ is as in the proof of Theorem 4.4.

## 5  Discrete semi Mittag-Leffler laws and Geometric Compounding

**Definition 5.1** : A r.v $X$ is said to have a discrete semi Mittag-Leffler (DSML($a,b,\alpha$)) distribution if its PGF is given by $P(s) = [1+\psi(1-s)]^{-1}$, where $\psi(1-s) = a\,\psi[b(1-s)]$ for all $0<s<1$ and some $0<b<1<a$ satisfying $ab^{\alpha}=1$ for a unique $0<\alpha\leq 1$. This is the discrete analogue of semi Mittag-Leffler laws discussed in Sandhya (1991b).

**Definition 5.2** : A r.v $X$ is said to have a discrete Mittag-Leffler (DML($\lambda,\alpha$)) distribution if its PGF is given by $P(s) = [1+\lambda(1-s)^{\alpha}]^{-1}$, $\lambda>0$, $0<\alpha\leq 1$.

In this section, by a geometric($p$) variable we mean a geometric variable with PGF $ps/(1-qs)$, $q = 1-p$.

**Theorem 5.1** : A discrete r.v $X$ is a geometric($p$) sum of its own D-type variables, if and only if it is DSML($\frac{1}{p},b,\alpha$).

**Proof** : If $P(s)$ is the PGF of $X$, a DSML variable, then

$$P(s) = \frac{1}{1 + \psi(1-s)}$$

$$= \frac{1}{1 + a\psi[b(1-s)]}$$

$$= \frac{p}{p + \psi[b(1-s)]}, \quad p = \tfrac{1}{a}$$

$$= \frac{p/[1 + \psi(b(1-s))]}{1 - q/[1 + \psi(b(1-s))]},$$

which proves our assertion in both the directions.

**Theorem 5.2** : If a discrete r.v is a geometric ($p$) sum of its own D-type for two values of $p$, say $p_1$ and $p_2$ such that $\log p_1 / \log p_2$ is irrational, then it is DML($\lambda, \alpha$).

**Proof** : $\log p_1 / \log p_2$ is irrational implies $\log a_1 / \log a_2$ is irrational where $a_i = 1/p_i$, $i=1,2$. Hence $\psi(1-s) = \lambda(1-s)^\alpha$ for $\lambda > 0$ (as in Theorem.4.3) and the variable is DML($\lambda, \alpha$).

**Note** : Sandhya and Satheesh (1996) showed that a semi-$\alpha$-Laplace law is in class-L if and only if it is $\alpha$-Laplace. Restricting the support to $[0, \infty)$ and then to the non-negative lattice $I_o$, we have the following two results invoking Corollary 2.2.

**Theorem 5.3** : A semi Mittag-Leffler law is in class-L if and only if it is Mittag-Leffler.

**Theorem 5.4** : A DSML law is in discrete class-L if and only if it is DML.

**Theorem 5.5** : Consider a sequence $\{Y_j\}$ of i.i.d Bernoulli variables with parameter $b = p^{1/\alpha}$, $0<p<1$, $0<\alpha \leq 1$. Let $M$ be a non-negative lattice variable independent of $\{Y_j\}$, with PGF P($s$) such that

$$\lim_{s \to 1} \frac{1 - P(s)}{(1-s)^\alpha} = \lambda > 0, 0 < \alpha \leq 1 \tag{1}$$

and put $X = \sum_{j=1}^{M} Y_j$. Let $\{X_i\}$ be a sequence of independent copies of $X$ and define $S_N = \sum_{i=1}^{N} X_i$, where $N$ is a geometric($p$) variable independent of $\{X_i\}$. Then as $p \to 0$, $S_N$ converges in law to a DML variable.

**Proof** : The PGF of $S_N$ is given by

$$P_p(s) = \frac{pP[1 - b(1 - s)]}{1 - qP[1 - b(1 - s)]}$$

$$= \frac{P[1 - b(1 - s)]}{P[1 - b(1 - s)] + p^{-1}\{1 - P[1 - b(1 - s)]\}}$$

Now, $\dfrac{1 - P[1 - b(1 - s)]}{p} = \dfrac{1 - P[1 - p^{1/\alpha}(1 - s)]}{[p^{1/\alpha}(1 - s)]^{\alpha}}(1 - s)^{\alpha}$.

Under the condition (5.1) the R.H.S converges to $\lambda(1-s)^{\alpha}$ as $p \to 0$. Hence

$$\lim_{p \to 0} P_p(s) = [1 + \lambda(1-s)^{\alpha}]^{-1},$$

proving the result.

**Theorem 5.6** : In the set up of Theorem 5.5, let $\{Y_j\}$ be Bernoulli with parameter $b$. Then $S_N$ is identically distributed as $M$ if and only if $b^{\alpha} = p$ for a unique $0 < \alpha \leq 1$ and $M$ is DML($\lambda, \alpha$).

**Proof** : We have $P(s) = \dfrac{pP[1 - b + bs]}{1 - qP[1 - b + bs]}$.

Setting $P(s) = \phi(1-s)$ and $\psi(1-s) = \dfrac{1}{\phi(1 - s)} - 1$, we have

$$\frac{1}{1 + \psi(1 - s)} = \frac{p/[1 + \psi(b(1 - s))]}{1 - q/[1 + \psi(b(1 - s))]}$$

$$= \frac{1}{1 + a\psi(b(1 - s))}, \quad a = \tfrac{1}{p}.$$

Hence $P(s)$ corresponds to a DSML($\frac{1}{p}$, $b$, $\alpha$), with $\alpha$ defined by $b^\alpha = p$ for a unique $0 < \alpha \leq 1$. Writing $Q(1-s) = \dfrac{\psi(1-s)}{(1-s)^\alpha}$, we have:

$$\frac{1}{1+(1-s)^\alpha Q(1-s)} = \frac{1}{1+ab^\alpha(1-s)^\alpha Q(b(1-s))}.$$

As $ab^\alpha = 1$, we have $Q(1-s) = Q(b(1-s))$. On iteration $Q(1-s) = Q(b^n(1-s))$ for each positive integer $n$. Since $0 < b < 1$, his means $Q(1-s) = \lambda > 0$ under the condition (5.1). Hence $P(s) = [1 + \lambda(1-s)^\alpha]^{-1}$.

**Note** : Notice that the geometric parameter $p$ and the Bernoulli parameter $b$ are related by $b^\alpha = p$ and the choice is under the restriction $b \leq p$ (c.f. Remark 4.1).

**Theorem 5.7** : In the setup of Theorem 5.6, let $\{M_i\}$ be a sequence of independent copies of $M$, and $N_0$ be a geometric($p_0$) variable independent of $M$ and $p_0 \neq p$. Then $\sum_{i=1}^{N_0} M_i$ and $\sum_{i=1}^{N} X_i$ are identically distributed if and only if $p < p_0$, $b^\alpha = p/p_0$ and $M$ is DML($\lambda, \alpha$).

**Proof** : The condition is equivalent to

$$p_0\, P(s) / (1 - q_0\, P(s)) = [p\, P(1-b+bs)] / [1-[q\, P(1-b+bs)]]$$

Setting $\phi(1-s)$, $\psi(1-s)$ and $Q(1-s)$ as in Theorem 5.6 we have

$$p_0^{-1}\, \psi(1-s) = p^{-1}\, \psi[b(1-s)].$$

Thus $0 < \alpha \leq 1$ is uniquely defined by $b^\alpha = p/p_0$. This shows that $p < p_0$ as $0 < b < 1$. Further,

$$p\,(1-s)^\alpha\, Q(1-s) = p_0\, b^\alpha (1-s)^\alpha\, Q[b(1-s)]$$

and hence $Q(1-s) = Q[b(1-s)]$. Now, proceeding as in the proof Theorem 5.6 we see that $P(s) = [1 + \lambda(1-s)^\alpha]^{-1}$ under the condition (5.1).

**Acknowledgement**

The authors wish to thank both the referees and the editor Professor M B Rajarshi for their suggestions that has brought more clarity to the paper.

**S. Satheesh and N. Unnikrishnan Nair**

Department of Statistics

Cochin University of Science and Technology

Cochin - 682 022. India.

ssatheesh@sancharnet.in

**Notes added.**

1. We wish to acknowledge that Jayakumar (1995) has priority to our Theorem.5.1 since his semi-$\alpha$-geometric laws is the same as our discrete semi Mittag-Leffler ($a$, $b$,$\alpha$) laws. Jayakumar, K (1995). The stationary solution of a first order integer valued autoregressive process, *STATISTICA*, **anno LV**, n.2, 221 – 228.

2. See also Satheesh, Nair and Sandhya (2002). Stability of random sums, *Stochastic modeling and Applications*, **5**, p.17-26; (available at arxiv as math.PR/0311348) for a related and more general discussion.